# Quantum cohomology of flag manifolds $G/B$ and quantum Toda lattices

By BUMSIG KIM*

## 1. Introduction

Let $G$ be a connected semi-simple complex Lie group, $B$ its Borel subgroup, $T$ a maximal complex torus contained in $B$, and $\text{Lie}(T)$ its Lie algebra. This setup gives rise to two constructions; the generalized nonperiodic Toda lattice ([28], [29]) and the flag manifold $G/B$.

The Toda lattice for $(G, B, T)$ is the dynamical system on the cotangent bundle $T^*\text{Lie}(T)$ endowed with the canonical holomorphic symplectic form and the holomorphic hamiltonian function we consider in this paper,

$$(1) \qquad H(p,q) = (p,p) - \sum_{\text{simple roots } \alpha_i} (\alpha_i, \alpha_i) \exp(\alpha_i(q)),$$

where $(\,,\,)$ is any fixed nonzero multiplication of the Killing form on each simple component of $\text{Lie}(G)$ and the simple roots are given by the roots of $B$ with respect to $T$. This system is known to be completely integrable ([30], [29]). Therefore the variety defined by the ideal generated by the integrals of motions is the lagrangian analytic submanifold of $T^*\text{Lie}(T)$.

On the other hand, for the flag manifold $G/B$ we have the small quantum cohomology ring $QH^*(G/B, \mathbb{C})$ which we study in this paper.

The quantum cohomology of a compact symplectic manifold $X$ is $H^*(X, \mathbb{C})$ with a multiplication depending on the parameter space $H^*(X, \mathbb{C})$ ([26], [11]). The rigorous establishments are finished due to tremendous works [13], [32], [33], [34], [36], [37], [38], [39], [41] for symplectic manifolds and [4], [5], [6], [14], [25], [26], [31] for smooth projective varieties. Let us restrict the parameter space to $H^2(X, \mathbb{C})$ and denote this quantum cohomology by $QH^*(X, \mathbb{C})$ which will be called the small quantum cohomology ([14]). In fact, the relevant parameter space for the small quantum cohomology is $H^2(X, \mathbb{C})/H^2(X, 2\pi\sqrt{-1}\mathbb{Z})$.

Since $QH^*(G/B, \mathbb{C})$ is generated by second cohomology classes and parameters, it is natural to ask for the relations between those classes in the quantum cohomology ring. We find the answer. Let $q_i$ be the coordinates

*Partially supported by a Sloan Doctoral Dissertation Fellowship and an Institut Mittag-Leffler Postdoctral Fellowship.

of the parameter space $H^2(G/B, \mathbb{C})$ defined by $q_i(\sum a_j p_j) = \exp(a_i)$, where $p_j$ is the cohomology class corresponding to the fundamental weights according to Borel's description ([7]). If $l$ is the rank of $G$, then $QH^*(G/B, \mathbb{C})$ is $H^*(G/B, \mathbb{C}) \otimes_{\mathbb{C}} \mathbb{C}[q_1, \ldots, q_l]$ with a deformed ring structure such that modulo $q$'s the quantum ring is exactly the ordinary ring. We obtain a ring *presentation* of $QH^*(G/B)$.

THEOREM I. *The small quantum cohomology ring $QH^*(G/B, \mathbb{C})$ is canonically isomorphic to*

$$\mathbb{C}[p_1, \ldots, p_l, q_1, \ldots, q_l]/I$$

*where $I$ is the ideal generated by the nonconstant complete integrals of motions of the Toda lattice for the Langlands-dual Lie group $(G^\vee, B^\vee, T^\vee)$ of $(G, B, T)$.*

Some explanations are in order. The Langlands-dual Lie group $G^\vee$ is the complex connected Lie group, that has the root system $(G^\vee, B^\vee, T^\vee)$ dual to $(G, B, T)$ and $\text{Lie}(T^\vee) = \text{Lie}(T)^*$ ([16]).

Let $\{T_k\}$ and $\{T^k\}$ be dual bases of $H^*(G/B)$ in the sense that $\int_{G/B} T_i T^j = \delta_{i,j}$. Since $G/B$ is a homogeneous space, the quantum product $x \circ y$ of two cohomology classes $x$ and $y$ can be written as

$$x \circ y = \sum_{\beta \in H_2(G/B), k} q^\beta T_k I_\beta(x, y, T^k),$$

where $I_\beta(x, y, T^k)$ is the number of discrete rational curves $C$ in $G/B$ with $[C] = \beta$ passing through three cycles Poincaré-dual to $x$, $y$ and $T^k$ and

$$q^\beta = \Pi_{i=1}^l q_i^{<p_i, \beta>}.$$

These $I_\beta(x, y, T^k)$, the so-called 3-pointed genus 0 Gromov-Witten invariants, encode the enumerative geometry of $G/B$ for rational curves. Since Theorem I is a ring presentation, it does not provide the invariants, just as the simple presentation of the classical cohomology of Grassmannians does not yield the multiplication table of the Schubert cells.

In fact we prove more in this paper. Givental considered quantum $\mathcal{D}$-modules, which play a significant role in the mirror theorem ([17], [18], [19]). *In particular, using the quantum hyperplane section principle ([2], [3], [24]) it is possible to compute the virtual numbers of rational curves in Calabi-Yau 3-fold complete intersections in homogeneous spaces with the knowledge of the quantum $\mathcal{D}$-module structure of the ambient spaces.* We show that the $\mathcal{D}$-module structure for $G/B$ is governed by the conservation laws of *quantum* Toda lattices which are the quantizations of the Toda lattices and still integrable ([27], [28], [35], [40]). The hamiltonian operator we consider is

$$\hat{H} = \Delta - \sum_{\text{simple roots } \alpha_i} (\alpha_i, \alpha_i) \exp(\alpha_i(q)),$$

where $\Delta$ is the Laplacian on $\mathrm{Lie}\,T$ associated with the invariant form $(\,,\,)$. Let $\mathcal{D}$ be the differential operator algebra over $\mathbb{C}$, generated by $\hbar\frac{\partial}{\partial t_i}$, multiplication by $\hbar$ and $\exp t_i$.

THEOREM II. *The quantum $\mathcal{D}$-module of $G/B$ is canonically isomorphic to $\mathcal{D}/\mathcal{I}$ where $\mathcal{I}$ is the left ideal generated by the nonconstant complete quantum integrals of motions of the quantum Toda lattice for the Langlands-dual Lie group $(G^\vee, B^\vee, T^\vee)$ of $(G, B, T)$.*

*Remark.* 0. Our result is true for $\mathbb{Q}$ in place of $\mathbb{C}$ as the coefficient ring of cohomology.

1. Using the Quot schemes for flags as a compactification of $\mathrm{Mor}(\mathbb{P}^1, \mathrm{SL}(n,\mathbb{C})/B)$, Ciocan-Fontanine ([8]) proved Theorem I in the case of $G = \mathrm{SL}(n,\mathbb{C})$, which was conjectured in [21].

2. Fomin-Gel'fand-Postnikov ([12]) obtained the structure constants of the small quantum cohomology for $\mathrm{SL}(n,\mathbb{C})/B$.

3. Givental ([20]) proved his mirror conjecture for $\mathrm{SL}(n,\mathbb{C})/B$ based on quantum Toda lattices.

4. In his unpublished works Dale Peterson was able to obtain all 3-pointed genus 0 Gromov-Witten invariants for any $G/P$ where $P$ is a parabolic subgroup. In particular, Theorem I is also a result of his works.

*Structure of the paper.* In Section 2, we recall the classical Toda lattices and the construction of integrals. In Section 3, we review the definition of quantum cohomology, quantum $\mathcal{D}$-modules and state two properties of quantum $\mathcal{D}$-modules which will be key ideas in proving theorems. In Section 4, we prove Theorem II and obtain Theorem I by a specialization of Theorem II. In the proof of Theorem II, the integrability of quantum Toda lattices will be crucial. Two examples for $G = \mathrm{SL}(n,\mathbb{C})$ and $G = \mathrm{SO}(2n+1,\mathbb{C})$ are given in Section 5. In Section 6, the equivariant quantum cohomology of the flag manifold is formulated by the lagrangian foliation of the Toda dynamical system on $T^*\mathrm{Lie}\,(T^\vee)$.

## 2. Toda lattices

In this section we recall the construction of integrals of motions of the Toda lattice for $(G, B, T)$ as defined in the introduction. We follow [28], [29]. To do so, we first describe the coadjoint method in symplectic geometry and factor the symplectic space $T^*\mathrm{Lie}\,T$ through one of the coadjoint orbits of $B$. By using Ad-invariant polynomial functions followed by a twist we obtain Poisson commutative functions including the Hamiltonian $H$ of (1) which we are interested in.

For $\chi \in (\operatorname{Lie} T)^*$ let $t_\chi \in \operatorname{Lie} T$ be the unique vector such that $\alpha(t_\chi) = (\alpha, \chi)$, $\alpha \in (\operatorname{Lie} T)^*$, let $E_\alpha$ denote a nonzero root vector corresponding to the root $\alpha$ (so that $[t, E_\alpha] = \alpha(t) E_\alpha$ for any $t \in \operatorname{Lie} T$), $l = \operatorname{rank} G$, and let $\pi_i$ be the weights such that $(\pi_i, \alpha_j) = \delta_{i,j}$. If we regard $e = \sum_{j=1}^l E_{-\alpha_j}$ as an element in $(\operatorname{Lie} B)^* = (\operatorname{Lie} T \oplus \sum_{\alpha > 0} \mathbb{C} E_\alpha)^*$ using the form $(,)$, the coadjoint orbit $B \cdot e$ of $e$ is seen to be $\operatorname{Lie} T \oplus \sum_{j=1}^l \mathbb{C}^\times E_{-\alpha_j}$ endowed with the induced holomorphic KKS symplectic form $\sum_{i=1}^l d[t_{\pi_i}] \wedge \frac{d[E_{\alpha_i}]}{[E_{\alpha_i}]}$ from the Poisson manifold $(\operatorname{Lie} B)^* \cong \operatorname{Lie} B' := \operatorname{Lie} T \oplus_{\alpha < 0} \mathbb{C} E_\alpha$. This orbit $B \cdot e$ is the augmented dynamical system of the cotangent space $T^* \operatorname{Lie} T$ in the Toda lattice under the map $\xi$:

$$
\begin{array}{rcl}
(\operatorname{Lie} T)^* \oplus \operatorname{Lie} T & \overset{\xi}{\to} & \operatorname{Lie} T \oplus \sum_{j=1}^l \mathbb{C}^\times E_{-\alpha_j} \\
(\lambda, \sum y_i t_{\pi_i}) & \mapsto & (t_\lambda, \sum e^{y_i} E_{-\alpha_i}).
\end{array}
\tag{2}
$$

The above map is a symplectic morphism. The hamiltonian function $H$ will be factored through a function on the coadjoint orbit.

For the integrals of motions (including the augmented $H$) consider the set $\operatorname{Poly}(\operatorname{Lie} G)^{\operatorname{Ad}}$ of all Ad-invariant polynomial functions (i.e., Casimir elements) $u$ on $\operatorname{Lie} G$, the restriction $u \in \operatorname{Poly}(\operatorname{Lie} B')$, and then $f$-twisting: $u^f(x) = u(x + f)$, for $x \in \operatorname{Lie} B'$, where

$$
f := -\frac{1}{2} \sum_i (\alpha_i, \alpha_i) E_{\alpha_i}.
\tag{3}
$$

Let us think of

$$
u_T^f := u^f \circ \xi \in \mathcal{O}((\operatorname{Lie} T)^* \oplus \operatorname{Lie} T).
\tag{4}
$$

Notice that $(,)_T^f$ augments $H$ of (1) if

$$
(E_{\alpha_i}, E_{-\alpha_i}) = 1,
\tag{5}
$$

which we will always assume. In the Toda lattice we thus obtain many integrals of motions which are pairwise commutative. We will see this well-known fact in Lemma 2. Since $\operatorname{Lie} G / \operatorname{Ad} = \operatorname{Lie} T / W$ (here $W$ denotes the Weyl group) is isomorphic to $\mathbb{C}^{\operatorname{rank} G}$ (due to Chevalley's theorem), we conclude that the dynamical system is completely integrable in the Liouville sense (see [30] for the proof of the linear independence of the differentials of some $l$-many integrals in involution with each other).

## 3. Quantum cohomology

The quantum cohomology multiplication is a certain formal Frobenius supermanifold $H^*(X, \mathbb{C})$ ([26], [11]). The manifold $M := H^*(X)$ has the vector space structure, a metric $g$ given by the Poincaré pairing, and the structure sheaf $\mathcal{O}_M$ of formal $\mathbb{C}$-valued series with respect to the vector space structure. The metric is just a nondegenerate symmetric $(0,2)$-tensor field. Notice that the metric is compatible to the vector space structure, i.e. $g(X,Y)$ is constant for constant fields $X$, $Y \in \Gamma(M, \mathcal{T}_M)$. To define the supercommutative multiplication $\circ : S^2(\mathcal{T}_M) \to \mathcal{T}_M$ we recall basic facts on Gromov-Witten theory, which are highly nontrivial.

A prestable rational curve is, by definition, a connected *arithmetic* genus 0 projective curve with, at worst, nodes as singular points. It could be reducible. Consider a prestable map $f$ from $C$ to $X$ with fixed ordered $n$-many marked points $x_i \in C$. We will identify $(f, C, \{x_i\})$ with $(f', C', \{x_i'\})$ if there is an isomorphism $h$ from $C$ to $C'$ preserving the configuration of marked points such that $f = f' \circ h$. A stable map $(f, C, \{x_i\})$ is a prestable map with only finitely many automorphisms. In other words, when $f$ maps an irreducible component of $C$ to a point in $X$, the component has at least three special points (nodes and marked points).

Let $X_{n,\beta}$ be the coarse moduli space of all stable maps to $X$ from prestable curves $C$ with $n$-ordered marked points with the fixed second cohomology $\beta = f_*([C]) \in H_2(X, \mathbb{Z})$. (See [14] for details.) Now assume that $X$ is a homogeneous projective variety for a complex reductive Lie group. Then it is convex since its holomorphic tangent space is generated by global holomorphic sections. According to [25], [6], [14], $X_{n,\beta}$ — whenever it is nonempty — is a compact complex *orbifold* with complex dimension $\dim X + <c_1(T_X), \beta> + n - 3$. More precisely, locally near a stable map the moduli space has data of a quotient of a holomorphic domain acted by the (finite) group of all automorphisms of the stable map. There are natural morphisms on the moduli spaces, namely, evaluation maps $e_i$ at $i^{\text{th}}$ marked points and forgetting-marked-point maps $\pi$:

$$\begin{array}{ccc} X_{n+1,\beta} & \xrightarrow{e_{n+1}} & X \\ \pi \downarrow & & \\ X_{n,\beta}, & & \end{array}$$

where $X_{n+1,\beta}$ is the universal family (not in the strict sense, but it is "with the ambiguity of the automorphisms") of $X_{n,\beta}$ whenever the moduli space is nonempty. If $s_i$ are the universal sections for the marked points, then $e_i = e_{n+1} \circ s_i$. Here we assume that $\pi$ is the forgetful map of the last marked point.

Let us define the potential $\Phi \in \mathcal{O}_M(H^*(X))$ (here $H^*(X)$ is considered as a formal scheme),

$$\Phi(\gamma) = \sum_{n \geq 3} \sum_{\beta \in H_2(X)} \frac{1}{n!} \int_{X_{n,\beta}} e_1^*(\gamma) \ldots e_n^*(\gamma),$$

which encodes all Gromov-Witten (GW) invariants for genus zero in one single generating function. The GW-invariant

$$\int_{X_{n,\beta}} e_1^*(A_1)\ldots e_n^*(A_n)$$

is equal to the number of stable maps from $\mathbb{P}^1$ to $X$ passing through general cycles Poincaré-dual to $A_i$ at $i^{\text{th}}$ marked points. Notice that the valid sum over $\beta$ in the potential is finite for each $n$ because $<c_1(T_X),\beta>$ is positive if $\beta$ is a nonzero effective class. $\Phi$ is a $\mathbb{C}$-valued formal series on $H^*(X)$ considered as a linear space. Now we are ready to define the quantum multiplication. For the constant vector fields $X,Y,Z$ the product $\circ$ is defined to satisfy $g(X \circ Y, Z) = XYZ(\Phi(\gamma))$. Here $Z(\Phi)$ is the derivative of the formal function $\Phi$ along the vector field $Z$. Since $\Phi(\gamma)$ (or its derivatives) may not converge in general, $X \circ Y$ is a formal section of the tangent sheaf $\mathcal{T}_M$. The condition of the associativity of the multiplication is exactly the WDVV-equation which is satisfied because of the boundary structure of the moduli space $X_{n,\beta}$. The quantum cohomology has the identity 1: $1 \circ X = X$.

To complete the Frobenius manifold picture of quantum cohomology we need an Euler vector field. Let $\{p_a, p_i\}$ be a linear basis of $H^*(X)$ and $\{\beta_i\} \subset H_2(X,\mathbb{C})$ be the dual basis of the second cohomology classes $\{p_i\}$. If $t_a$ denotes the coordinates with respect to the basis of the cohomology space, an Euler vector field is defined to be $E = \sum_a (1-\deg(p_a))t_a \frac{\partial}{\partial t_a} + \sum_i <c_1(T_X), \beta_i> \frac{\partial}{\partial t_i}$, and so $E(\Phi) = (3-\dim X)\Phi$, modulo up to quadratic terms. It can be seen that $L_E g = (2-\dim X)g$ and $L_E \circ = \circ$.

In the Frobenius manifold $M$, the quantum multiplication $\circ$ defines a connection:

$$\nabla_A(B) = \hbar \nabla_A^0(B) - A \circ B,$$

where $\nabla^0$ is the Levi-Civita connection of $g$ and $\hbar$ is a formal parameter. It is flat (for each $\hbar$) ([11]). (Strictly speaking, $\frac{1}{\hbar}\nabla$ is a connection.)

Now let us restrict the parameter space $H^*(X)$ to $H^2(X)$. This means the multiplication is defined by the following simple form: For $A, B, C \in H^*(X)$ and $p \in H^2(X)$

$$g(A \circ_p B, C) = \sum_{\beta \in H_2(X)} e^{<p,\beta>} \int_{X_{3,\beta}} e_1^*(A) \wedge e_2^*(B) \wedge e_3^*(C).$$

The reason is the so-called divisor axiom, i.e., $\pi_* e_{n+1}^*(p) = <p,\beta>$.

Now assume that the basis $\{p_i\}$ is chosen in closed Kähler cone. Let $q_i$ be functions on $H^2(X)$ defined by $q_i(\sum a_j p_j) = \exp a_i$.

*Definition.* Since $A \circ B \in H^*(X,\mathbb{C}) \otimes_\mathbb{C} \mathbb{C}[q_1,\ldots,q_l]$, the *small* quantum cohomology $QH^*(X)$ is defined to be $H^*(X) \otimes_\mathbb{C} \mathbb{C}[q_1,\ldots,q_l]$ with $\circ$ multiplication.

Notice that due to the Euler vector field $QH^*(X)$ is a graded algebra, where the cohomology classes are given their usual complex degrees and $q_j$ has degree $<c_1(T_X),\beta_j>$.

The kernel of the natural homomorphism from $\operatorname{Sym}(H^2(X))\otimes_\mathbb{C}\mathbb{C}[q_1,\ldots,q_l]$ to $QH^*(X)$ yields the *characteristic subscheme* of $T^*H^2(X)$ ([21], [1]). The map is surjective. One of the goals of the paper is to describe the characteristic subscheme.

Next we recall Givental's quantum $\mathcal{D}$-module, of which we will make use in order to prove our claim in the introduction. If

$$s \in \Gamma(H^2(X), \mathcal{O}_{H^*(X)})(= H^*(X)\otimes\mathbb{C}[[t_1,\ldots,t_l]]$$

with the chosen basis $\{p_i\}\subset H^2(X)$ and its dual $\{\beta_j\}\subset H_2(X)$, $l=\operatorname{rank} H_2(X)$), the quantum differential equations of parallel sections are

$$\hbar\frac{\partial}{\partial t_i}s = p_i \circ s, \ i=1,\ldots,l,$$

whose fundamental solutions can be constructed ([9], [11], [19]). Let $c$ be the first Chern classes of the line bundle over $X_{n,\beta}$ where fibers are tangent lines at the first marked point of prestable curves. This so-called universal line bundle at $i^{\text{th}}$ marked points is the normal bundle of the embedding of $X_{n,\beta}$ under the section $s_i$. For $A\in H^*(X)$ define an $H^*(X)$-valued function $s_A$ depending on formal variables $t_i$, and $\hbar^{-1}$ using the Poincaré pairing with $B\in H^*(X)$ as follows:

$$g(s_A, B) = g(A\exp pt/\hbar, B) + \sum_{d\neq 0}\exp dt \int_{X_{2,d\beta}} \frac{e_1^*(A)\exp(e_1^*(pt)/\hbar)}{\hbar+c}\wedge e_2^*(B),$$

where we use multi-indices, $dt = \sum d_i t_i$, $d\beta = \sum d_i \beta_i$, $e_1^*(pt) = \sum e_1^*(p_i)t_i$. One can show that they form fundamental solutions. Now we want to focus on

$$S_A = g(s_A, 1) \in \mathbb{C}[t_1,\ldots,t_l,\hbar^{-1}][[q_1=\exp t_1,\ldots,q_l=\exp t_l]].$$

The quantum $\mathcal{D}$-module of $X$ is defined to be the $\mathcal{D}$-module generated by these formal functions $S_A$, $A\in H^*(X)$, where $\mathcal{D}$ is the operator algebra over $\mathbb{C}$, generated by $\hbar\frac{\partial}{\partial t_i}$, multiplication by $\hbar$ and $\exp t_i$.

*Definition.* The quantum $\mathcal{D}$-module of $X$ is

$$\mathcal{D}/\mathcal{I} \text{ where } \mathcal{I} = \{D\in\mathcal{D}: DS_A = 0 \text{ for all } A\in H^*(X)\}.$$

Let us give degrees

$$\deg\frac{\partial}{\partial t_i} = 0, \ \deg\hbar = 1, \ \deg\exp t_i = <c_1(T_X),\beta_i>, \ \deg t_i = 0.$$

Then $\mathcal{I}$ *is homogeneous* since $S_A$ is homogeneous if $A$ is a homogeneous class.

Suppose a differential operator $P(\hbar\frac{\partial}{\partial t_i}, \exp t_i, \hbar)$ generated by $\hbar, \hbar\frac{\partial}{\partial t_i}, \exp t_i$ annihilates $S_A$ for all $A \in H^*(X)$. Applying $P$ to $S_A = g(s_A, 1)$ one concludes that $P(p_i \circ, q_i, 0) = 0$ in the small quantum cohomology (see Cor. 6.4 in [19]). Conversely, if one has $P(p_i \circ, q_i, 0) = 0$ which is a quadratic polynomial in each $p_i$, then the differential operator $P(\hbar\frac{\partial}{\partial t_i}, \exp t_i, 0)$ ordering $\exp t_i$ before the positions of $\hbar\frac{\partial}{\partial t_i}$ annihilates $S_A$ for all $A$. *These two facts ([19], [20]) will be used in the next section to prove our claims in the introduction.*

## 4. Quantum $\mathcal{D}$-module of $G/B$

In this section $X = G/B$ will be considered. This space is a complex homogeneous projective variety. We may find an embedding into the complex projective space of the finite dimensional irreducible highest weight representation space $V(\rho)$ where $\rho$ is the half sum of all positive roots. The embedding space is the orbit of the highest weight vector in $\mathbb{P}(V(\rho))$.

The classical cohomology algebra $H^*(G/B)$ has Borel's description (see, for instance, [7]). It is generated by the second cohomology classes with relations generated by Weyl group $W$ invariant nonconstant classes. Notice that $G/B = K/(B \cap K)$, where $K$ is a maximal compact subgroup of $G$ containing $T$, and $W := N(T \cap K)/(T \cap K)$. So there is a $W$-action (induced by the right action of $N(T \cap K)$) on the cohomology group.

We find the quantum-counterpart of these relations. These relations on $QH^*(G/B)$ are provided by integrals of motions of the Toda lattices for the Langlands-dual Lie group $G^\vee$ to $G$. The Langlands-dual group $G^\vee$ is the connected complex Lie group that has $T^\vee := \mathrm{Hom}(\pi_1(T), \mathbb{C}^\times)$ as a maximal complex torus and has the root system of $(G^\vee, T^\vee)$ dual to that of $(G, T)$. Canonically $\mathrm{Lie}(T^\vee) \cong (\mathrm{Lie}\, T)^*$. (See [16] for a reference.) Recall that $(\mathrm{Lie}\, T^\vee)^* \ni \alpha_i^\vee$ are roots dual to $\alpha_i \in (\mathrm{Lie}\, T)^*$, and $\alpha_i^\vee := \frac{2(\alpha_i, \cdot)}{(\alpha_i, \alpha_i)}$. For example, the root system of simple Lie algebra type $B_n$ is dual to that of $C_n$ and all other types are self-dual. The author uses the Langlands-dual group rather than the dual Lie algebra in order to emphasize the canonical way of the identifying $\mathrm{Lie}\, T^\vee$ with $(\mathrm{Lie}\, T)^*$ and to indicate the plausible relevance of loop groups and representation theory to mirror symmetry phenomena.

To simplify the notation, let $\mathrm{Lie}\, G = \mathfrak{g}$, $\mathrm{Lie}\, T = \mathfrak{h}$, $\mathrm{Lie}\, B = \mathfrak{b}$. We first construct the ring homomorphism $\phi$:

$$
\begin{aligned}
\mathrm{Poly}((\mathfrak{h}^\vee)^*) \otimes \mathcal{O}(\mathfrak{h}^\vee) &\supset \mathrm{Sym}(\mathfrak{h}^*) \otimes \mathbb{C}[e^{\alpha_1^\vee}, \ldots, e^{\alpha_l^\vee}] \\
&\xrightarrow{\phi} QH^*(G/B) \\
&= H^*(G/B) \otimes_\mathbb{C} \mathbb{C}[q_1, \ldots, q_l],
\end{aligned}
\tag{6}
$$

and show that the kernel of $\phi$ contains the integrals of the Toda lattice on $\mathfrak{h}^{\vee*} \oplus \mathfrak{h}^\vee$ with the Hamiltonian $H^\vee(p, q) = (p, p) - \sum(\alpha_i^\vee, \alpha_i^\vee) \exp(\alpha_i^\vee(q))$.

We begin with Borel's description of the ordinary cohomology ring $H^*(G/B,\mathbb{C})$. Given an integral weight $\chi$, there is a character $\theta \in \mathrm{Mor}(T,\mathbb{C}^\times)$ such that $\theta(\exp t) = \exp \chi(t)$ for $t$ in the Lie algebra $\mathrm{Lie}\,T$ of $T$; extend $\theta$ to a character of $B$ by setting $\theta(n) = 1$ for $n$ in the unipotent radical $N$ of $B$. This is possible since $B$ is the semidirect product of $N$ and $T$. Then $G \times_B \mathbb{C} = G \times \mathbb{C}/\sim \,\to G/B$ is a complex line bundle where $(g,x) \sim (gb, \theta(b)x)$ for $g \in G$, $x \in \mathbb{C}$, $b \in B$. Sending the integral weight $\chi$ to the first Chern class $c_\chi$ of the complex line bundle, we obtain a ring homomorphism from the symmetric algebra $\mathrm{Sym}(\mathfrak{h}^*)$ to $H^*(G/B,\mathbb{C})$. This homomorphism is known to be surjective with the kernel generated by Weyl group $W$ action invariant elements without constant terms. By the choice of the Weyl chamber to make the roots of $B$ positive, the Chern classes $c_{\lambda_i}$ associated to the fundamental weights $\lambda_i$ are then in the closed Kähler cone; therefore $\{c_{\lambda_i}\}$ forms a $\mathbb{Z}_+$-basis of the cone. This depends on our choice of the definition of $G \times_B \mathbb{C}$. Now it is obvious how to define $\phi$: just identify the integral weight $\chi$ with the Chern class $c_\chi$. Then the basis in closed Kähler cone is $\{p_i = c_{\lambda_i}\}$. One sees that $\phi$ is onto.

From $(G,B,T)$ we naturally obtain $(G^\vee, B^\vee, T^\vee)$ and so on for other dualizations of the notation. If $\mathfrak{n}$ denotes $[\mathfrak{b},\mathfrak{b}]$, we may reconsider $u_T^f$ of (4) as an element of $\mathrm{Sym}(\mathfrak{b}^\vee/[\mathfrak{n}^\vee,\mathfrak{n}^\vee])$ using the form $(\,,\,)$. One sees that from $\xi$ of (2) and $\phi$ of (2) $\varphi := \phi \circ \xi^*$:

$$\mathrm{Sym}(\mathfrak{b}^\vee/[\mathfrak{n}^\vee,\mathfrak{n}^\vee]) = \mathbb{C}[[\lambda_1],\dots,[\lambda_l],[E_{\alpha_1^\vee}],\dots,[E_{\alpha_l^\vee}]] \to QH^*(G/B,\mathbb{C})$$

is the $\mathbb{C}$-algebra homomorphism mapping $[E_{\alpha_i^\vee}]$ and $[\lambda_i]$ to $q_i$ and the first Chern class $c_{\lambda_i}$ of the line bundle associated with the $i^{\mathrm{th}}$ fundamental weight $\lambda_i$ of the Lie group $G$, respectively. This is nicely summarized in the diagram:

$$\begin{array}{ccccc}
\mathrm{Sym}(\mathfrak{b}^\vee/[n^\vee,n^\vee]) & \xrightarrow{\xi^*} & \mathrm{Poly}((\mathfrak{h}^\vee)^*) \otimes \mathbb{C}[e^{\alpha_j^\vee}] & \xrightarrow{\phi} & QH^*(G/B) \\
[\lambda_i],\ [E_{\alpha_j^\vee}] & \longrightarrow & \lambda_i,\ e^{\alpha_j^\vee} & \longrightarrow & p_i = c_{\lambda_i},\ q_j.
\end{array}$$

Notice that

$$(\,,\,)_T^f = \sum_{1 \le i,j \le l}(\alpha_i^\vee,\alpha_j^\vee)[\lambda_i][\lambda_j] - \sum_{i=1}^l (\alpha_i^\vee,\alpha_i^\vee)[E_{\alpha_i^\vee}].$$

PROPOSITION.  $(\,,\,)_T^f$ *is in the kernel of* $\varphi$.

*Proof.* First, in the usual cohomology ring $\sum(\alpha_i^\vee,\alpha_j^\vee)c_{\lambda_i}c_{\lambda_j} = 0$. In the small quantum cohomology $c_{\lambda_i} \circ c_{\lambda_j} = c_{\lambda_i}c_{\lambda_j} + \delta_{i,j}q_i$ since $c_{\lambda_i}$ are dual to $\beta_j$, $\{c_{\lambda_i}\}$ is $\mathbb{Z}_+$-basis in the closed Kähler cone of $H^2(X,\mathbb{Z})$, and the fibers of the quotient map $G/B \to G/P_{\alpha_i}$ are all the rational curves of type $\beta_i$, where $P_{\alpha_i}$ is the parabolic subgroup containing $B$ and having only $-\alpha_i$ for negative roots. □

If $\{p_i\}$ is the basis $\{c_{\lambda_i}\}$ of $H^2(X, \mathbb{Z})$, the fundamental solutions $S_A$ of the quantum differential equations can be constructed as in Section 3. Let $H_\partial$ be the Laplace operator $\sum_{i,j}(\alpha_i^\vee, \alpha_j^\vee)\hbar^2 \frac{\partial}{\partial t_i}\frac{\partial}{\partial t_j}$, and $H_q$ the differential operator $-\sum_i (\alpha_i^\vee, \alpha_i^\vee)\exp t_i$ (just multiplication), and define

$$\hat{H} := H_\partial + H_q.$$

Then $\hat{H}S_A = 0$ because of the observation at the end of Section 3 and the above proposition.

Put $a_i = (\overbrace{0, \ldots, 1}^{i}, \ldots, 0) \in \mathbb{Z}^l$. Following [20], we have a recursive formula of $S_A$.

COROLLARY. $S_A = g(e^{pt/\hbar}s, A)$, where $s = \sum s^{(d)}\exp(dt)$ and $s^{(d)} \in H^*(X)$, has the recursion relation

$$s^0 = 1, \ \hbar[\hbar\sum_{i,j}d_i d_j(\alpha_i^\vee, \alpha_j^\vee) + 2\sum_{i,j}d_i(\alpha_i^\vee, \alpha_j^\vee)p_j]s^{(d)} = \sum_{i\ |d_i|>0}(\alpha_i^\vee, \alpha_i^\vee)s^{(d-a_i)}.$$

Even though the differential operators are not commutative, we use the notation, say, $\mathbb{C}[P, Q]$, for the $\mathbb{C}$-algebra generated by operators $P$, $Q$.

LEMMA 1. *Suppose $\hat{H}s = 0$ for a formal function $s = \sum_{d_i \geq 0} f_d(t, \hbar)\exp dt$, where $f_d(t, \hbar)$ is a polynomial in $t = (t_1, \ldots, t_l)$, $\hbar$, and $\hbar^{-1}$. Let $D \in \mathbb{C}[\hbar\frac{\partial}{\partial t_i}, \exp t_i, \hbar, \hbar^{-1}]$ be another operator such that $[\hat{H}, D] = 0$ and $Ds \equiv 0$ modulo $(\exp t_1, \exp t_2, \ldots, \exp t_l)$. Then $Ds = 0$.*

*Proof.* Since $\hat{H}Ds = 0$, it suffices to show that if $\hat{H}s = 0$ and $s \equiv 0$ modulo $(\exp t_1, \exp t_2, \ldots, \exp t_l)$, then $s = 0$. We may let $s = \sum_{|d|=k>0}f_d(t, h)\exp dt +$ higher order terms in $\exp t$, where $|d| := d_1 + \cdots + d_l$. Pick any $d$ such that $|d| = k$. Since $\hat{H} = H_\partial + H_q$, $H_\partial f_d(t)\exp dt = 0$. But $H_\partial f_d(t, h)\exp dt = f_d(t, h)H_\partial \exp dt + f'_d(t, h)\exp dt + f''_d(t, h)\exp dt$, where $f'_d$ and $f''_d$ have lower degrees than $f_d(t, h)$ in $t$. Hence $f_d(t, h)H_\partial \exp dt = 0$. But since $(,)$ is a nonzero multiple of a positive definite metric in the $\mathbb{R}$-space $\mathfrak{h}_\mathbb{R}^*$ spanned by all roots — $(,)$ is chosen to be a nonzero multiple of the Cartan-Killing form (we may assume that $G$ is a simple Lie group and then a product of simple Lie groups) —, $H_\partial \exp dt$ is nonzero unless $d = (0, \ldots, 0)$. Therefore $f_d = 0$ for $|d| = k$. This implies $s = 0$. □

For the sake of simple notation, $S(\mathfrak{g})$ will denote the symmetric algebra $\text{Sym}(\mathfrak{g})$ of $\mathfrak{g}$. Let $\eta : S(\mathfrak{g}^\vee)^{G^\vee} \to S(\mathfrak{h}^\vee)^W$ be the quotient map by the ideal generated by $n^\vee$ and $(n')^\vee$. We need an elaborate version of the integrability of quantum Toda lattices, namely

LEMMA 2 ([27], [35], [40]). *There is an algebra injective homomorphism $\hat{\gamma}_f$ from $S(\mathfrak{g}^\vee)^{G^\vee}$ to $\mathbb{C}[\hbar\frac{\partial}{\partial t_1},\ldots,\hbar\frac{\partial}{\partial t_l},\exp t_1,\ldots,\exp t_l,\hbar]$ satisfying the following properties*:

0) *If $u \in S_k(\mathfrak{g}^\vee)^{G^\vee}$, then $\hat{\gamma}_f(u)$ is degree $k$ homogeneous with given degrees 0, 1, 2, to $\frac{\partial}{\partial t_i}$, $\hbar$, $\exp t_i$.*

1) *The principal symbol of $\hat{\gamma}_f(u)$, i.e. replacing $\hbar\frac{\partial}{\partial t_i}$, $\exp t_i$, $\hbar$ with $[\lambda_i]$, $[E_{\alpha_i^\vee}]$, 0, respectively, is $u_T^f$.*

2) *When $u$ is the invariant form $(\,,\,) \in S_2(\mathfrak{g}^\vee)^{G^\vee}$, $\hat{\gamma}_f(u)$ is $\hat{H}$. In particular, $[\hat{\gamma}_f(u), \hat{H}] = 0$ for any $u \in S(\mathfrak{g}^\vee)^{G^\vee}$.*

3) *For any $u \in S(\mathfrak{g}^\vee)^{G^\vee}$, $\hat{\gamma}_f(u)$ modulo $(\exp t_1,\ldots,\exp t_l)$ is $\eta(u)$ if we replace $\hbar\frac{\partial}{\partial t_i}$ with $\lambda_i$.*

The proof of this lemma will be given later following [27], [35], [40]. Denote $\hat{\gamma}_f(u)$ by $\hat{u}_T^f$. Now we are fully ready to prove

COROLLARY. *The generators $S_A$ of the quantum $\mathcal{D}$-module of $G/B$ vanish under the nonconstant operators $\hat{u}_T^f$.*

*Proof.* Consider $u_T^f$ as in Lemma 2, and notice that $(\hat{u}_T^f)S_A$ vanishes modulo $(\exp t_1, \exp t_2,\ldots,\exp t_l)$ because of 3) of Lemma 2. We can apply lemma 1 with $s = S_A$ and $D = \hat{u}_T^f$ since $[\hat{u}_T^f, \hat{H}] = 0$ as in 2) of Lemma 2. We conclude that $(\hat{u}_T^f)S_A = 0$. □

Result 1) of Lemma 2 yields Theorem I. (The injectivity of the induced surjective morphism $\mathbb{C}[p,q]/I$ to $QH^*(G/B)$ can be shown using Poincaré series of graded $\mathbb{C}$-modules.)

THEOREM II. *The quantum $\mathcal{D}$-module of $G/B$ is canonically isomorphic to the quotient of $\mathcal{D}$ by the left ideal $\mathcal{I}$ generated by $\hat{\gamma}_f(S_k(\mathfrak{g}^\vee)^{G^\vee})$, $k > 0$.*

*Proof.* Suppose that $D \in \mathcal{D}$ and $DS_A = 0$ for all $A \in H^*(X)$. We will show that $D$ is in $\mathcal{I}$. First, we may assume that $D$ has a homogeneous degree. If the degree of $D$ is 0, it is done. Use the induction on the degree of $D$. We may assume that $D$ is not an element of $\hbar\mathcal{D}$, otherwise it is done. Since the symbol of $D$ is a relation in the small quantum cohomology, the symbol is $vu_T^f$ for some $u \in S(\mathfrak{g}^\vee)^{G^\vee}$ and some $v \in \mathbb{C}[p_1,\ldots,p_l,,q_1,\ldots,q_l]$. Then $D' := D - \hat{v}\hat{u}_T^f$ is in $\hbar\mathcal{D}$, where $\hat{v}$ is any differential operator with the symbol $v$. Since $D'S_A = 0$ for all $A \in H^*(X)$, by induction we deduce that $D'$ is in $\mathcal{I}$ and hence $D$ is in $\mathcal{I}$. □

*Proof of Lemma* 2. For the sake of simplicity we will prove this lemma for $\mathfrak{g}$ instead of $\mathfrak{g}^\vee$. Recall that we denote by $\mathfrak{h}$ a Cartan subalgebra, by $\alpha_i$ simple roots, and by $W$ the Weyl group. Let $\Phi^+$ be the set of all positive roots, and $\mathfrak{b} := \mathfrak{h} \oplus \sum_{\beta \in \Phi^+} \mathbb{C} E_\beta$, $\mathfrak{b}' := \mathfrak{h} \oplus \sum_{\beta \in \Phi^+} \mathbb{C} E_{-\beta}$. Notice that if $(\pi_i, \alpha_j) = \delta_{i,j}$, then $[t_{\pi_i}, E_{\alpha_j}] = \alpha_j(t_{\pi_i}) E_{\alpha_j} = \delta_{i,j} E_{\alpha_j}$. Thus if $\mathfrak{n} = [\mathfrak{b}, \mathfrak{b}]$, the quotient Lie algebra $\mathfrak{b}/[\mathfrak{n}, \mathfrak{n}]$ is isomorphic to a Lie algebra generated by $P_i$ and $Q_j$ with relations $[P_i, Q_j] = \delta_{i,j} Q_j$, $[P_i, P_j] = 0 = [Q_i, Q_j]$.

*Step* 1. First let us look into the universal enveloping algebra $U(\mathfrak{b})$ and then $U(\mathfrak{b}/[\mathfrak{n}, \mathfrak{n}])$. Using the invariant form $(,)$ which is in each simple component a nonzero multiple of the Cartan-Killing form, we will identify $\mathfrak{b}^*$ (resp. $\mathfrak{n}^*$) with $\mathfrak{b}'$ (resp. $\mathfrak{n}' = [\mathfrak{b}', \mathfrak{b}']$). Consider an algebra homomorphism $f : U(\mathfrak{n}') \to \mathbb{C}$, and denote by $U_f$ the kernel of $f$. For the moment, $f$ is not necessarily $f$ in the statement of Lemma 2. In fact, $-\frac{1}{2} \sum_i (\alpha_i, \alpha_i) E_{\alpha_i}$ defines a character of $\mathfrak{n}'$ using the form $(,)$, and hence a homomorphism $f : U(\mathfrak{n}') \to \mathbb{C}$ which will be used at the end of the proof. Notice that $U := U(\mathfrak{g}) = U(\mathfrak{b}) \otimes U(\mathfrak{n}') = U(\mathfrak{b}) \otimes (\mathbb{C}1 \oplus U_f) = U(\mathfrak{b}) \oplus UU_f$. This decomposition allows us to have a $\mathbb{C}$-linear projection $\delta_f$ from $U(\mathfrak{g})$ to $U(\mathfrak{b})$ along $UU_f$. The linear map $\delta_f$ restricted to the center $Z(\mathfrak{g})(=: Z)$ is a ring homomorphism. Denote by $\gamma_f$ the map $\delta_f$ followed by the $\rho$-twisting. Here the $\rho$-twisting is an automorphism of $U(\mathfrak{b})$ defined by $X \to X + \frac{1}{2} \mathrm{tr}(\mathrm{ad}_{\mathfrak{b}} X)$, $X \in \mathfrak{b}$. We then have the following facts on the ring homomorphism $\gamma_f : Z \to U(\mathfrak{b})$:

a) The image of $\gamma_0$ (i.e. when $f = 0$) is in $U(\mathfrak{h})$, (this is the so-called Harish-Chandra isomorphism in literature ([10])).

b) The composition $\sigma \circ \gamma_f$ is nothing but $\gamma_0$, where $\sigma$ is the "modulo-$\mathfrak{n}$" map from $U(\mathfrak{b})$ to $U(\mathfrak{h})$ along $\mathfrak{n}U(\mathfrak{b})$ (notice that $\mathfrak{n}U(\mathfrak{b})$ is a two sided ideal of $U(\mathfrak{b})$),

c) The image of $\gamma_0$ is in $S(\mathfrak{h})^W$. Furthermore, $\gamma_0 : Z \to S(\mathfrak{h})^W$ is an isomorphism.

Let $x_0 \in \mathfrak{h}$ be the unique element given by $<\alpha_i, x_0> = 1$ for all simple roots $\alpha_i$. By the extended derivation $\mathrm{ad}\, x_0$ of $S(\mathfrak{g})$ and $U(\mathfrak{g})$, Kazhdan ([27], [35]) introduced the $\mathrm{ad}\, x_0$-grading (resp. $\mathrm{ad}\, x_0$-filtration) in $S(\mathfrak{b}) = \bigoplus_{k=0}^\infty S_{(k)}(\mathfrak{b})$ (resp. $U(\mathfrak{b}) = \sum_{k=0}^\infty U_{(k)}(\mathfrak{b})$). By definition $S_{(i)}(\mathfrak{b})_j$ is the eigenspace of $\mathrm{ad}\, x_0|_{S_i(\mathfrak{b})}$ corresponding to the eigenvalue $j$, and $S_{(k)}(\mathfrak{b}) = \bigoplus_{i+j=k} (S_{(i)}(\mathfrak{b}))_j$. The $\mathrm{ad}\, x_0$-filtration of $U(\mathfrak{b})$ is obtained in a similar way. Notice that $f$ is determined by values on $\mathfrak{n}'/[\mathfrak{n}', \mathfrak{n}']$, and hence $f$ defines $\sum_i a_i E_{\alpha_i}$, $a_i \in \mathbb{C}$. For $u \in S_k(\mathfrak{g})^G$ let us define $u^f \in S(\mathfrak{b}) = \mathrm{Poly}(\mathfrak{b}')$ by $u^f(x) = u(x + \sum_i a_i E_{\alpha_i})$, $x \in \mathfrak{b}'$. Then in addition to a), b), c), following Kostant ([27]) we claim:

d) i): $\gamma_f(\gamma_0^{-1} \circ \eta(u)) \in U_{(k)}(\mathfrak{b})$ and ii): if $f$ is nonsingular, i.e. $f(E_{-\alpha_i}) \neq 0$ for all simple roots $\alpha_i$, $\gamma_f$ is an injection of $\mathrm{ad}\, x_0$-filtered algebras such that for $u \in S_k(\mathfrak{g})^G$, $\tau_{(k)}(\gamma_f \circ \gamma_0^{-1} \circ \eta(u)) = u^f$, where $\tau_{(k)}$ is the unique map from

$U_{(k)}(\mathfrak{b})$ to $S_{(k)}(\mathfrak{b})$, $k = 0, 1, \ldots$, such that the sequence

$$0 \longrightarrow U_{(i-1)}(\mathfrak{b}) \xrightarrow{\text{injection}} U_{(i)}(\mathfrak{b}) \xrightarrow{\tau_{(i)}} S_{(i)}(\mathfrak{b}) \longrightarrow 0$$

is exact and if $v \in U_{(i)}(\mathfrak{b})$, $w \in U_{(j)}(\mathfrak{b})$, then $\tau_{(i+j)}(vw) = \tau_{(i)}(v)\tau_{(j)}(w)$.

Here are proofs of the claims.

*Proof of* a). Let $v \in U$. Since $U(\mathfrak{b}) = U(\mathfrak{h}) \oplus \mathfrak{n}U(\mathfrak{b})$, it is enough to show that the component of $v$ in $\mathfrak{n}U$ is also in $U\mathfrak{n}' = UU_0$. According to the PBW theorem, $X_\alpha^q H^m X_{-\alpha}^p$ form a basis of the enveloping algebra $U$, where we assume that $p, q, m$ are multi-indices for the nonzero multi-vectors $X_{(-)\alpha}$ corresponding to the positive (negative) roots $(-)\alpha$, and the multi-vectors $H$ form a basis of Cartan subalgebra. If $v = X_\alpha^q H^m X_{-\alpha}^p$, then $[t, v] = 0$ for $t$ in the Cartan subalgebra will imply $(q - p)\alpha(t) = 0$ in each basis. Hence if $v$ is any element in $Z \cap \mathfrak{n}U$, then it is also in $U\mathfrak{n}'$.

*Proof of* b). (Prop. 2.3. [27]). Let $v \in Z$ and $v_0 = \delta_0(v)$. Then $v - v_0 \in U\mathfrak{n}'$ and also in $\mathfrak{n}U$, because $Z \cap U\mathfrak{n}' = Z \cap \mathfrak{n}U$ seen in the proof of a). Thus $v - v_0$ is in the first component of the decomposition $U = (\mathfrak{n}U + UU_f) \oplus U(\mathfrak{h}) = (\mathfrak{n}U(\mathfrak{b}) \oplus UU_f) \oplus U(\mathfrak{h})$ (the reason of the first decomposition is $\mathfrak{n}U \subset \mathfrak{n}U(\mathfrak{b}) + UU_f$ and the second one which is seen before). Since $v - v_0 \in nU(\mathfrak{b}) \oplus UU_f$, $\delta_f(v) - \delta_0(v) \in nU(\mathfrak{b})$. Hence $\gamma_f(v) - \gamma_0(v) \in nU(\mathfrak{b})$.

*Proof of* c). Let $\rho$ be half the sum of all the positive roots and let $v$ be an element of the center $Z$ of the universal enveloping algebra $U$. It can be uniquely written as $v = v_0 + v_1$, $v_0 \in U(\mathfrak{h})$ and $v_1 \in UU_0 = U\mathfrak{n}'$. Considering $U(\mathfrak{h}) = S(\mathfrak{h}) = Poly(\mathfrak{h}^*)$, we need to show only $v_0(\lambda + \rho) = v_0(s_\alpha \lambda + \rho)$, where $-\lambda \in \mathfrak{h}^*$ is a dominant weight, $\alpha$ is a simple root of $\mathfrak{g}$, and $s_\alpha$ is the reflection associated with $\alpha$. Consider the Verma module $M(\lambda) = U(\mathfrak{g}) \otimes_{U(\mathfrak{b}')} \mathbb{C}$. Here $U(\mathfrak{g})$ is a right $U(\mathfrak{b}')$-module and $\mathbb{C}$ is a twisted-$\mathfrak{b}'$-module by $(b_\mathfrak{h} + b_{\mathfrak{n}'})x = (\lambda + \rho)(b_\mathfrak{h})x$, where $b_\mathfrak{h} \in \mathfrak{h}$, $b_{\mathfrak{n}'} \in \mathfrak{n}'$, $x \in \mathbb{C}$. Let $X$ be the canonical generator of $M(\lambda)$ and let $Y$ be the canonical generator of $M(s_\alpha \lambda)$. Then $v_0(\lambda + \rho)X = v_0 X = (v_0 + v_1)X = vX$, and $v_0(s_\alpha \lambda + \rho)Y = vY$. But it is known that $M(s_\alpha \lambda)$ is a sub-$\mathfrak{g}$-module of $M(\lambda)$, and the action of the center $Z$ on $M(\lambda)$ is "scalar multiplications." Hence we conclude the proof. The proof that $\gamma_0$ is an isomorphism onto $S(\mathfrak{h})^W$ can be found in, for instance, [10] and [42].

Notice that if $v$ is the universal Casimir element

$$\sum_i t_{\alpha_i} t_{\pi_i} + \sum_{\alpha \in \Phi^+} \frac{E_\alpha E_{-\alpha}}{(E_\alpha, E_{-\alpha})} + \sum_{\alpha \in \Phi^+} \frac{E_{-\alpha} E_\alpha}{(E_\alpha, E_{-\alpha})},$$

then

(7) $$\gamma_0(v) = \sum_i t_{\alpha_i} t_{\pi_i} - (\rho, \rho),$$

which is in $S(\mathfrak{h})^W$.

*Proof of* d). The claim without $\rho$-twisting is given in Kostant's paper ([27]) (see Theorems 2.4.1 and 2.4.2). Since the $\rho$-twisting preserves the filtrations and $\tau_{(i)} \circ \rho-\text{twisting} = \tau_{(i)}$, d) is true.

From now on, let $f = -\frac{1}{2}\sum_i (\alpha_i, \alpha_i) E_{\alpha_i}$, which is nonsingular.
Claims a), b), c) and d) show that:

e) For $u \in S_k(\mathfrak{g})^G$, $\gamma_f \circ \gamma_0^{-1} \circ \eta(u) \in U_{(k)}(\mathfrak{b})$ are involutions and $\tau_{(k)}(\gamma_f \circ \gamma_0^{-1} \circ \eta(u)) = u^f$, and $\gamma_f \circ \gamma_0^{-1} \circ \eta(u)$ modulo $\mathfrak{n}$ is $\eta(u)$.

*Step* 2. Define a new Lie algebra $\mathfrak{b}_\hbar := \mathfrak{b} \otimes \mathbb{C}[\hbar]$, $[x \otimes \hbar^i, y \otimes \hbar^j]_\hbar := [x,y]\hbar^{i+j+1}$ for $x$, $y$ in $\mathfrak{b}$, and give the standard *grading* on $U(\mathfrak{b}_\hbar)$. Since the ideal generated by $x \otimes y - y \otimes x - [x,y]_\hbar$ is homogeneous in the graded tensor $\mathbb{C}[\hbar]$-algebra $T(\mathfrak{b}_\hbar) := \bigoplus_{i=0}^\infty (\sum_{i+j=k} \mathfrak{b}^i \otimes \hbar^j)$, there is no problem to give the induced grading on $U(\mathfrak{b}_\hbar) := \bigoplus_{k=0}^\infty U_k(\mathfrak{b}_\hbar)$. We introduce a new *grading* on $U(\mathfrak{b}_\hbar)$ using an extension of $\text{ad}\, x_0$ to $U(\mathfrak{b}_\hbar)$ given by $\text{ad}\, x_0(\hbar^k x_1 \ldots x_m) = \hbar^k \sum_{i=1}^m x_1 \ldots \text{ad}\, x_0(x_i) \ldots x_m$ for $x_i \in \mathfrak{b}$. $U(\mathfrak{b}_\hbar)$ has a *grading* in a similar way: $U(\mathfrak{b}_\hbar) = \bigoplus_{k=0}^\infty U_{(k)}(\mathfrak{b}_\hbar)$ where $U_{(k)}(\mathfrak{b}_\hbar) = \bigoplus_{i+j=k}(U_i(\mathfrak{b}_\hbar))_j$, $(U_i(\mathfrak{b}_\hbar))_j$ is the eigenspace of $\text{ad}\, x_0|_{U_i(\mathfrak{b}_\hbar)}$ corresponding to the eigenvalue $j$. The linear map from $\mathfrak{b}_\hbar$ to $U(\mathfrak{b})$ sending $x \otimes \hbar^k$ to $x$ provides a unique homomorphism of algebras

$$\phi : U(\mathfrak{b}_\hbar) \to U(\mathfrak{b})$$

with $\phi(x \otimes \hbar^k) = x$. Notice that $\phi$ restricted to $U_{(k)}(\mathfrak{b}_\hbar)$ is an isomorphism onto $U_{(k)}(\mathfrak{b})$ because one can construct the inverse map. In fact, if $U_{(k)}^i(\mathfrak{b}_\hbar)$ is the quotient image in $U_{(k)}(\mathfrak{b}_\hbar)$ of $\bigoplus_{j \geq k-i} T(\mathfrak{b}) \otimes \hbar^j$, then there is a linear filtration $U_{(k)}(\mathfrak{b}_\hbar) = \sum_{i=0}^k U_{(k)}^i(\mathfrak{b}_\hbar)$ such that $\text{gr}\,_i U_{(k)}(\mathfrak{b}_\hbar) := U_{(k)}^i(\mathfrak{b}_\hbar)/U_{(k)}^{i-1}(\mathfrak{b}_\hbar) = S_{(i)}(\mathfrak{b})$.
For $u \in S_k(\mathfrak{g})^G$, let

$$\hat{u}^f = \phi^{-1}(\gamma_f \circ \gamma_0^{-1} \circ \eta(u)) \in U_{(k)}(\mathfrak{b}_\hbar).$$

Then they are involutions and $\tau_{(k)}(\phi(\hat{u}^f)) = u^f$. Notice also that $\hat{u}^f$ modulo $\mathfrak{n}$ is $\eta(u)$.

From the quotient maps $\mathfrak{b} \to \mathfrak{b}/[\mathfrak{n},\mathfrak{n}]$ and $\mathfrak{b}_\hbar \to (\mathfrak{b}/[\mathfrak{n},\mathfrak{n}])_\hbar$, we have an obvious commutative diagram:

$$\begin{array}{ccccc}
U_{(k)}(\mathfrak{b}_\hbar) & \xrightarrow{\phi} & U_{(k)}(\mathfrak{b}) & \xrightarrow{\tau_{(k)}} & S_{(k)}(\mathfrak{b}) \\
\downarrow & & \downarrow & & \downarrow \\
U_{(k)}((\mathfrak{b}/[\mathfrak{n},\mathfrak{n}])_\hbar) & \xrightarrow{\phi_T} & U_{(k)}(\mathfrak{b}/[\mathfrak{n},\mathfrak{n}]) & \xrightarrow{\tau_{(k),T}} & S_{(k)}(\mathfrak{b}/[\mathfrak{n},\mathfrak{n}]).
\end{array}$$

We give the analogous ad $x_0$-filtration/grading on $U((\mathfrak{b}/[\mathfrak{n},\mathfrak{n}])_{(\hbar)})$ and $S(\mathfrak{b}/[\mathfrak{n},\mathfrak{n}])$. Consider $u_T^f$ and $\hat{u}_T^f$, which are the images of $u^f$ and $\hat{u}^f$ under the vertical maps in the above diagram, respectively. Thus finally we obtain $\hat{\gamma}_f(u) := \hat{u}_T^f \in U_{(k)}((\mathfrak{b}/[\mathfrak{n},\mathfrak{n}])_\hbar)$ such that $\tau_{(k),T} \circ \phi(\hat{u}_T^f) = u_T^f$, $\hat{u}_T^f$ commutes with each other. Additionally, $\hat{u}_T^f$ modulo $\mathfrak{n}/[\mathfrak{n},\mathfrak{n}]$ is $\eta(u)$. Regarding $U_{(k)}((\mathfrak{b}/[\mathfrak{n},\mathfrak{n}])_\hbar)$ as the differential operator algebra generated by $\hbar\frac{\partial}{\partial t_i}$ ($\leftrightarrow t_{\pi_i}$), $\exp t_j$ ($\leftrightarrow E_{\alpha_j}$), $\hbar(\leftrightarrow \hbar)$, we can state that $\tau_{(k),T} \circ \phi$ (resp. modulo $\mathfrak{n}/[\mathfrak{n},\mathfrak{n}]$-map) corresponds to the principal symbol map (resp. modulo $(\exp t_i)$-map). A simple check with (7) shows that

$$\hat{H} = \sum_i t_{\alpha_i} t_{\pi_i} - \sum (\alpha_i, \alpha_i) E_{\alpha_i} = \hat{(,)}_T^f.$$

□

## 5. Examples

We describe the formulation of our result for the complex simple Lie algebras of types $A$ and $B$. First some additional notation. For any basis $\{\beta_1, \ldots, \beta_l\}$ of $\mathfrak{h}^*$, let $\{p_i\}$ denote the corresponding basis of $H^2(G/B)$ and let $X(p,q) = \sum_i p_i t'_{\beta_i} + \sum_j q_j E_{-\alpha_j^\vee} + f$, where $(t_{\beta_i}, t'_{\beta_j}) = \delta_{i,j}$. It is required that $f = -\frac{1}{2}\sum_i (\alpha_i, \alpha_i) E_{\alpha_i}$ and $(E_{\alpha_i}, E_{-\alpha_i}) = 1$ as in (3) and (5). Denote by $E_{i,j}$ the matrix whose $(i,j)^{\text{th}}$ entry is 1 and all of whose other entries are 0. Write $\text{diag}(x_1, \ldots, x_n)$ for the diagonal matrix whose $(i,i)^{\text{th}}$ entry is $x_i$. If the symplectic Lie algebra $\text{sp}(2n, \mathbb{C})$ is defined to be the set of all complex $2n \times 2n$ matrices $Y$ such that $Y^T J + JY = 0$, where

$$J = \begin{pmatrix} 0 & 1 \\ -1 & 0 \end{pmatrix},$$

then $\text{sp}(2n, \mathbb{C})$ is the set of all matrices

$$(A, B, C) := \begin{pmatrix} A & B \\ C & -A^T \end{pmatrix},$$

where $B$ and $C$ are symmetric.

1) Suppose $G = \text{SL}(n+1, \mathbb{C})$, $B =$ the set of all upper triangle matrices with determinant 1. A Cartan subalgebra $\mathfrak{h}$ contained in $B$ is the set of all diagonal matrix $(x_1, \ldots, x_{n+1})$ such that $\sum x_i = 0$. Take weights $\varepsilon_i$ given by $\varepsilon_i(\text{diag}(x_1, \ldots, x_{n+1})) = x_i$. The simple roots $\alpha_i$ are thus $\varepsilon_i - \varepsilon_{i+1}$ in $n+1$ dimensional euclidian space $\mathbb{R}^{n+1}$ with the standard orthonormal basis $\varepsilon_i$. Therefore $\alpha_i^\vee = \alpha_i$ so that we may use the same Lie algebra for the dual Lie algebra, and $t_{\alpha_i} = E_{i,i} - E_{i+1,i+1}$, $E_{\alpha_i} = E_{i,i+1}$, $E_{-\alpha_i} = E_{i+1,i}$ ($\text{tr}(t_{\alpha_i} t_{\alpha_i}) = 2\text{tr}(E_{\alpha_i} E_{-\alpha_i})$). We conclude that $X(p,q) = \sum_i p_i t_{\alpha_i} + \sum_i q_i E_{-\alpha_i} - \sum_i E_{\alpha_i}$ is

$$\begin{pmatrix} p_1 & -1 & 0 & & 0 \\ q_1 & p_2 - p_1 & -1 & & \\ 0 & q_2 & \cdots & & 0 \\ & \cdots & \cdots & \cdots & -1 \\ 0 & & 0 & q_n & -p_n \end{pmatrix},$$

where $p_i$ stand for the cohomology classes associated to fundamental weights, and $q_j$ stand for the homology classes dual to $p_j$.

Consider the characteristic polynomials

$$\det(t + X) = t^{n+1} + \sum_{1 \le v \le n} J_v(X) t^{n-v}$$

of $X \in \mathrm{sl}(n+1, \mathbb{C})$. It is known that $\mathrm{Poly}(\mathfrak{g})^G = \mathbb{C}[J_1, \ldots, J_n]$. Hence we obtain the formula that the small quantum cohomology of a complete flag manifold $\mathrm{SL}(n+1, \mathbb{C})/B$ is isomorphic to the polynomial algebra $\mathbb{C}[p_i, q_i]$ quotient by the ideal generated by $J_v(X(p,q))$, $v = 1, \ldots, n$.

2) Let $G = \mathrm{SO}(2n+1, \mathbb{C})$, $n \ge 2$: The corresponding root space for the Lie algebra can be realized in the usual euclidian space $\mathbb{R}^n$ with the standard orthonormal basis $\varepsilon_i$ and simple roots $\alpha_i = \varepsilon_i - \varepsilon_{i+1}$ ($1 \le i \le n-1$), $\alpha_n = \varepsilon_n$. Let us use the symplectic Lie algebra $\mathrm{sp}(2n, \mathbb{C})$ for the dual Lie algebra with $\varepsilon_i(\mathrm{diag}(x_1, \ldots, x_n, -x_1, \ldots, -x_n)) = x_i$. We may then set $t_{\varepsilon_i} = (E_{i,i}, 0, 0)$, $E_{(-)\alpha_i^\vee} = (E_{i,i+1}, 0, 0)^{(T)}$, $i = 1, \ldots, n-1$, $E_{(-)\alpha_n^\vee} = (2)(0, E_{n,n}, 0)^{(T)}$. Therefore $X(x, q) = \sum_i x_i t_{\varepsilon_i} + \sum_i q_i (E_{i+1,i}, 0, 0) + 2 q_n (0, 0, E_{n,n}) - \sum_i (E_{i,i+1}, 0, 0) - 2(0, E_{n,n}, 0) =$

$$\begin{pmatrix} x_1 & -1 & 0 & & | & 0 & & & \\ q_1 & x_2 & -1 & & | & & & & \\ 0 & q_2 & \cdots & & | & & & & 0 \\ & & & x_n & | & & & 0 & -2 \\ - & - & - & - & + & - & - & - & - \\ 0 & & & & | & -x_1 & -q_1 & & \\ & & & & | & 1 & -x_2 & -q_2 & \\ & & & 0 & | & & 1 & \cdots & \\ & & 0 & 2q_n & | & & & & -x_n \end{pmatrix}.$$

It is known that the polynomials $J_v$ defined by

$$\det(t + X) = t^{2n} + \sum_{1 \le v \le n} J_v(X) t^{2(n-v)}$$

freely generate the invariant symmetric algebra $\mathrm{Poly}(\mathrm{sp}(2n, \mathbb{C}))^{\mathrm{Sp}(2n,\mathbb{C})}$. So the small quantum cohomology of the set of all complete orthogonal isotropy flags $\mathbb{C}^1 \subset \mathbb{C}^2 \subset \cdots \subset \mathbb{C}^n \subset \mathbb{C}^{2n+1}$ in $\mathbb{C}^{2n+1}$ is isomorhic to the polynomial algebra $\mathbb{C}[x_1, \ldots, x_n, q_1, \ldots, q_n]$ quotient by the ideal generated by $J_v(X(x,q))$, $v = 1, \ldots, n$. Here $x_i$ are the first Chern classes of the dual of the line bundles with fibers $\mathbb{C}^i/\mathbb{C}^{i-1}$ (recall our convention of $G \times_B \mathbb{C}$), and $q_j$ stand for the homology classes dual to $\sum_{i=1}^{j} x_j$, $j = 1, \ldots, n-1$, $\frac{x_1 + \cdots + x_n}{2}$.

## 6. Equivariant quantum cohomology

In this section $G$ denotes a connected *compact* Lie group. Let $X$ be a homogeneous complex algebraic manifold with a (hamiltonian) left $G$-action. The $G$-space $X$ gives rise to its homotopic quotient $X \times_G EG$, $(pg, q) \sim (p, gq)$ which is a total space of $X$-fiber bundle over $BG$. It is the universal $X$-bundle for the $G$-space $X$. Since the Leray spectral sequence of the universal $X$-bundle is degenerate ([15]), the equivariant Poincaré pairing $g$ is perfect (i.e., the induced map from $H^*_G(X)$ to $\mathrm{Hom}(H^*_G(X), H^*_G(pt))$ is an isomorphism), and the sequence

$$0 \to IH^*_G(X) \to H^*_G(X) \to H^*(X) \to 0$$

is exact, where $I$ is the augmented cohomology group of $H^*(BG)$. Of course, $H^*_G(X)$ is $H^*(BG)$-free. The equivariant cohomology of $G$-space $X$ is $H^*_G(X)$ endowed with the $H^*(BG)$-algebra structure.

First we consider equivariant GW-invariants to define the equivariant quantum cohomology ([19], [21], [22], [23]). The natural maps arising in the picture of GW-invariants are $G$-equivariant maps (e.g., sections, evaluation maps, forgetful maps, and contractions are equivariant):

$$\begin{array}{ccc} X_{n,\beta} \times_G EG & \xrightarrow{e_i} & X_G \\ \downarrow & & \\ \overline{\mathcal{M}}_n \times BG. & & \end{array}$$

One may think of $X^G_{n,\beta} := X_{n,\beta} \times_G EG$ as the space of all vertical stable maps to $X_G$ of degree $\beta$. The author would like to point out one of the properties which the equivariant GW-invariants have, the so-called divisor axiom [26]. If $p \in H^2_G(X)$, then one can show that $\pi_* e^*_{n+1} p = <i^*_X(p), \beta>$ by considering the following diagram from the inclusions:

$$\begin{array}{ccc} X & \xrightarrow{i_X} & X_G \\ \uparrow & & \uparrow \\ X_{n+1,\beta} & \longrightarrow & X^G_{n+1,\beta} \\ \downarrow & & \downarrow \pi \\ X_{n,\beta} & \longrightarrow & X^G_{n,\beta}, \end{array}$$

where $n \geq 1$ and $\beta \neq 0$. (If $i_X^*(p)$ is not a second cohomology class, put $<i_X^*(p),\beta>=0$.) For the validity of the equivariant WDVV-equation, see [19].

For $\gamma \in H_G^*(X)$ we can define an $H_G^*(pt)$-valued potential $\Phi(\gamma)$ by

$$\Phi(\gamma) = \sum_{n,\beta} \frac{1}{n!} <\gamma,\ldots,\gamma>_{n,\beta}^G,$$

where $<\gamma,\ldots,\gamma>_{n,\beta}^G$ is the *equivariant push forward* to $BG$ of the class $e_1(\gamma) \wedge \cdots \wedge e_n(\gamma)$ in $H_G^*(X_{n,G})$. Consider this as a formal power series with respect to an equivariant $H_G^*(pt,\mathbb{C})$-basis of $H_G^*(X)$. Define the equivariant quantum multiplication $\circ$ by the requirement $g(X \circ Y, Z) = XYZ\Phi$, $X, Y, Z \in H_G^*(X)$. The Euler vector field could be $\sum_a (1 - \deg(p_a)) t_a \frac{\partial}{\partial t_a} + \sum_i <c_1(T_X), \beta_i> \frac{\partial}{\partial t_i} + \sum \hbar_j \frac{\partial}{\partial \hbar_j}$. The variables $\hbar_j$ are from $G$-characteristic classes. The definition of the small equivariant quantum cohomology is obvious: this is the analogue of the small quantum cohomology.

*Example.* Let $S^1$ act on $\mathbb{P}^1$ as $[x:y] \mapsto [x \exp(2\pi\sqrt{-1}t) : y]$. Denote by $X_0$ and $X_\infty$ the equivariant Thom classes of the inclusions of the two fixed points $[1:0]$ and $[0:1]$, respectively. The small equivariant quantum cohomology of $\mathbb{P}^1$ is then

$$\mathbb{C}[X_0, X_\infty]/(X_0 X_\infty = q),$$

which follows from the grading of equivariant quantum cohomology and the usual relation $X_0 X_\infty = 0$ in the equivariant cohomology $H_{S^1}^*(\mathbb{P}^1)$.

Since $G/T \times_G EG = BT$, we have $H_G^*(G/T) = H^*(BT) = \mathrm{Sym}((\mathrm{Lie}\,T)^*)$ with $H^*(BG) = \mathbb{C}[c_1,\ldots,c_l] = H^*(BT)^W$ acting by $c_i = i^{\mathrm{th}}$ primary element. Denote by $\Sigma_i(t,q)$ the associated integrals of motions of the dual Toda lattice. The equivariant counterpart picture [19] for the quantum differential operators and its fundamental solutions leads to

THEOREM. *Let $G$ be a connected compact Lie group, $T$ a maximal torus, and $t_i$ the $G$-equivariant Euler classes of the line bundles over $G/T$ associated to fundamental weights with respect to a Weyl chamber. The $T$-equivariant small quantum cohomology of $G_\mathbb{C}/B$, where $B$ is the Borel subgroup of the complexification $G_\mathbb{C}$ of $G$, containing $T$, having the same Weyl chambers, is the quotient algebra $\mathbb{C}[t_1,\ldots,t_l,t_1',\ldots,t_l',q_1,\ldots,q_l]$ factored by relations,*

$$\Sigma_i(t,q) - c_i(t'),\ i=1,\ldots,l,$$

*where $q_i \in H_2(G/B, \mathbb{C}/2\pi\sqrt{-1}\mathbb{Z})$ are dual to the classes induced from $t_i$.*

Since $QH^*_G(X) \cong QH^*_T(X)^W$, one can obtain the above result also for the whole $G$-action. The reader might first derive $G$-equivariant one and then $T$-equivariant one using the restriction formula $QH^*_T(X) \cong QH^*_G(X) \otimes_{H^*_G(pt)} H^*_T(pt)$ ([23]). It seems that all operations on usual equivariant cohomology concerning only $G$ could be carried also in the quantum cases.

*Acknowledgment.* The author would like to thank E. Frenkel, H. Ooguri, N. Reshetikhin, and especially A. Givental for helpful discussions and the referees for useful suggestions.

University of California at Davis, Davis, CA
*E-mail address*: bumsig@math.ucdavis.edu